\documentclass[journal,transmag]{IEEEtran}
\usepackage{caption}
\usepackage{subcaption}
\usepackage{longtable}
\usepackage[sort, numbers]{natbib}
\usepackage[acronym]{glossaries}
\usepackage{graphicx}
\usepackage{amsmath}
\usepackage{amssymb}
\usepackage{breqn}
\usepackage{array}
\usepackage{xcolor}
\usepackage[capitalize]{cleveref}

\newacronym{abk:amod}{AMoD}{Autonomous Mobility-on-Demand}
\newacronym{abk:mpc}{MPC}{Model Predictive Control}
\newacronym{abk:mdp}{MDP}{Markov Decision Process}
\newacronym{abk:ro}{RO}{Robust Optimization}
\newacronym{abk:dro}{DRO}{Distributionally Robust Optimization}
\newacronym{abk:duro}{DURO}{Deep Uncertainty Robust Optimization}
\newacronym{abk:dnn}{DNN}{Deep Neural Networks}
\newacronym{abk:cnn}{CNN}{Convolutional Neural Network}
\newacronym{abk:lstm}{LSTM}{Long-Short-Term memory}
\newacronym{abk:gnn}{GNN}{Graph Neural Networks}
\newacronym{abk:gcn}{GCN}{Graph Convolution Network}
\newacronym{abk:gat}{GAT}{Graph Attention Network}
\newacronym{abk:so}{SO}{Stochastic Optimization}
\newacronym{abk:rl}{RL}{Reinforcement Learning}
\newacronym{abk:dohv}{DOHV}{Deterministic optimization with historical average demand}
\newacronym{abk:donn}{DONN}{Deterministic optimization with point prediction by a neural network}

\begin{document}

\title{Robust Vehicle Rebalancing with Deep Uncertainty in Autonomous Mobility-on-Demand Systems}

\author{{Xinling Li$^{1}$, Xiaotong Guo$^{1}$, Qingyi Wang$^1$, Gioele Zardini$^{1}$, and Jinhua Zhao$^2$}
\thanks{$^1$Department of Civil and Environmental Engineering, Massachusetts Institute of Technology, Cambridge, MA 02139 USA
        {\tt \{xinli831,xtguo,qingyiw,gzardini\}@mit.edu}}%
\thanks{$^2$Department of Urban Studies and Planning, Massachusetts Institute of Technology, Cambridge, MA 02139 USA
        {\tt jinhua@mit.edu}}%
\thanks{This work is supported by the US Department of Energy/Office of Energy Efficiency and Renewable Energy, Grant Agreement DE-EE0011186. This article
solely reflects the opinions and conclusions of its authors and not the Department of Energy, Office of Energy Efficiency and Renewable Energy,
or any other entity.}}

\maketitle

\begin{abstract}
\gls{abk:amod} services offer an opportunity for improving passenger service while reducing pollution and energy consumption through effective vehicle coordination. A primary challenge in the autonomous fleets coordination is to tackle the inherent issue of supply-demand imbalance. A key strategy in resolving this is vehicle rebalancing, strategically directing idle vehicles to areas with anticipated future demand. Traditional research focuses on deterministic optimization using specific demand forecasts, but the unpredictable nature of demand calls for methods that can manage this uncertainty. This paper introduces the \gls{abk:duro}, a framework specifically designed for vehicle rebalancing in \gls{abk:amod} systems amidst uncertain demand based on neural networks for robust optimization. \gls{abk:duro} forecasts demand uncertainty intervals using a deep neural network, which are then integrated into a robust optimization model. We assess \gls{abk:duro} against various established models, including deterministic optimization with refined demand forecasts and \gls{abk:dro}. Based on real-world data from New York City (NYC), our findings show that \gls{abk:duro} surpasses traditional deterministic models in accuracy and is on par with \gls{abk:dro}, but with superior computational efficiency. The \gls{abk:duro} framework is a promising approach for vehicle rebalancing in \gls{abk:amod} systems that is proven to be effective in managing demand uncertainty, competitive in performance, and more computationally efficient than other optimization models.
\end{abstract}


\section{Introduction} \label{introduction}

\IEEEPARstart{T}{here} has been a remarkable surge in the adoption of on-demand services such as Uber and Lyft during recent years. This growth is driven by factors such as heavy traffic and high cost of parking, which have diminished the attractiveness of private vehicle ownership. 
Industry reports project a compound annual growth rate exceeding 8.75\%~\citep{LIU2022100075}, estimating that the global on-demand market will reach an impressive value of \$218 billion by 2025 \citep{RN5}. Such services free users from the financial and logistical burdens of vehicle maintenance and parking, while offering the convenience of personal mobility.
The development of the paradigm of robo-taxi has been instrumental in creating platforms that centralize control and improve service quality. The fully compliant nature of autonomous vehicles opens the opportunity for effective fleet management that can mitigate issues such as pollution and traffic congestion \cite{tresca2025robo, kuhnimhof2023mobility}.

While the flexibility of \gls{abk:amod} services makes them highly appealing, it also introduces operational challenges~\citep{ZardiniAnnRev2022}.
Such services allow users to request rides from any location to their desired destination, with vehicles dispatched to meet the requests.
This passenger-centric approach adds layers of complexity to the operational planning process. 
Specifically, imbalances in vehicle availability frequently occur due to the uneven distribution of pick-up and drop-off locations. 
Addressing these challenges requires operational strategies that balance long-term system sustainability with user satisfaction and operator profitability.

To address the supply and demand imbalance in \gls{abk:amod} services, a common method involves rearranging the distribution of vehicles to align with the anticipated future demand. This process, known as \emph{vehicle rebalancing}, involves a strategic repositioning of vehicles over time to ensure that their availability is aligned with the predicted patterns of user trips. 
A key aspect of successful rebalancing strategies is the ability to forecast future demand accurately. 
Therefore, enhancing the precision of demand predictions or considering potential uncertainties in demand while making rebalancing decisions is essential.

The study of the demand forecasting model has been a significant topic of research for a long time. 
Much of the research on demand prediction has concentrated on achieving better accuracy in point forecasts~\citep{yao2018deep,du2020traffic,huang2022dynamical}. 
However, the unpredictable nature of real-world events means that future demand inherently involves uncertainty. 
Relying solely on precise point predictions for rebalancing optimization in \gls{abk:amod} services may not yield the best results due to this uncertainty, and the solutions derived may lack robustness. Therefore, it is crucial to factor in demand uncertainty during the optimization process to overcome these challenges.

The challenge of vehicle rebalancing in on-demand systems, especially considering demand uncertainty, has been a focus of various research methodologies. 
Existing approaches include modeling the problem as a \gls{abk:mdp}~\citep{ALKANJ20201088}, leveraging \gls{abk:mpc}~\citep{RAMEZANI2018203}, or data-driven models~\citep{9341481}. 
In \gls{abk:mpc} and data-driven methods, the optimization models are typically deterministic. To incorporate uncertainty, \gls{abk:ro} has been used to solve the rebalancing problem~\citep{GUO2021161}. 
\gls{abk:ro}~\citep{ben2009robust}, a significant discipline in operations research, is designed to address optimization problems with uncertainty. However, a limitation of \gls{abk:ro} is that it involves manual selection of the uncertainty sets. 
In contrast, \gls{abk:dro}~\citep{Rahimian_2022}, utilizing ambiguity sets derived by data-driven methods, offers a more customized approach. 
However, this method suffers from computational inefficiency and practical challenges, making it less suitable for real-time applications such as vehicle rebalancing.

In this study, we introduce a framework designed to improve rebalancing strategies in \gls{abk:amod} systems under demand uncertainty. 
Unlike traditional \gls{abk:ro} methods~\citep{GUO2021161}, which rely on manual input for defining uncertainty sets, our proposed framework leverages advanced deep neural networks to determine these intervals automatically. 
Moreover, it offers computational advantages over traditional data-driven \gls{abk:ro} methods.

To the best of our knowledge, this is among the first to integrate deep learning with \gls{abk:ro} methodologies to address demand uncertainty in \gls{abk:amod} operations. 
The contributions of this paper are threefold.
First, we introduce \gls{abk:duro}, a flexible framework suitable for decision-making problems under uncertainty.
By incorporating proposed deep uncertainty quantification techniques, \gls{abk:duro} enables more effective selection of uncertainty sets within the RO framework.
Second, we propose a matching-integrated rebalancing model to the vehicle rebalancing problem in \gls{abk:amod} operations, accounting for the uncertainty of the demand.
Finally, we benchmark \gls{abk:duro} against traditional point-prediction-based optimization, RO models, and DRO models.
Our findings demonstrate \gls{abk:duro}'s performance in ensuring system efficiency and robustness.

The remainder of this paper is organized as follows: Section \ref{literature} covers a literature review of the previous research in this field. 
\cref{method} introduces and formally presents the proposed \gls{abk:duro} framework. In \cref{experiment}, we present the experiment details and analyze the obtained results. \cref{conclusion} summarizes the main discoveries and provides an outlook on future research.


\section{Literature review} \label{literature}
Vehicle rebalancing has long been a focal point in the operations of on-demand systems. Numerous studies have been dedicated to developing more effective and practical methods to enhance operational quality. Typically, this research is bifurcated into two main areas: demand prediction and devising rebalancing strategies. In this section, we provide an overview of the literature in both these domains and pinpoint existing research gaps.

\subsection{Demand prediction}\label{demand}
Demand prediction based on advanced neural network structures and training algorithms is a dominant direction in this field. In transportation systems, the demand exhibits both spatial and temporal dimensions. Thus, it is essential to concurrently consider these dual aspects to attain a more precise prediction. To account for these features, a considerable body of research focused on improving demand prediction by the \gls{abk:cnn} \citep{lecun1998gradient} and \gls{abk:lstm} \citep{10.1162/neco.1997.9.8.1735}. The convolution structure of 
\gls{abk:cnn} can capture the spatial features of the demand data, and the \gls{abk:lstm} is used for modeling the temporal correlation\citep{yao2018deep}. 

One limitation of \gls{abk:cnn}s is that it is primarily designed for unstructured data like images, and it may not generalize effectively to structured data, such as trip demand data. To capture the spatial correlations in unstructured graphs, \gls{abk:gnn}~\citep{kipf2016semi} is developed as a generalization of \gls{abk:cnn} which can be adapt to handling complex graph structures. \gls{abk:gnn} offers enhanced customization in terms of the number of nodes and the format of the adjacency matrix, making them more suited for predicting \gls{abk:amod} demand than \gls{abk:cnn}s. Due to the aforementioned advantage, there are researchers focusing on adapting \gls{abk:gnn} for traffic prediction. For example, \cite{Yu_2018} applied \gls{abk:gcn} and \cite{zheng2019gman} applied \gls{abk:gat}, which are two widely used types of \gls{abk:gnn}, to traffic forecasting. They both demonstrated the superiority of \gls{abk:gnn} to other state-of-the-art models for demand prediction. To get a comprehensive introduction to the development of the point prediction models in traffic demand forecasting, interested readers can refer to the review papers by \cite{ghalehkhondabi2019review} and \cite{doi:10.1080/01441647.2018.1442887}.


While a considerable body of research in demand prediction has focused on improving the accuracy of point prediction, uncertainty persists in transportation systems and many other real-world systems, such as robotics and power grids. In particular, demand prediction constantly involves uncertainty, either due to the model error or data randomness. 
The oversight of uncertainty in demand prediction can lead to cumulative mispredictions in subsequent models~\citep{zhao2002propagation}. Therefore, it is essential to quantify the uncertainty during demand prediction. \cite{zhuang2022uncertainty, wang2023uncertainty} made probabilistic assumptions on the demand distribution and use a neural network to predict the parameters for the distribution. However, a single distribution might not be sufficient to fit the demand. \cite{rodrigues2020beyond} proposed to solve this problem by using quantile prediction to capture the uncertainty of traffic demand. By directly predicting quantiles, no distributional assumption needs to be made but the number of predicted quantiles is fixed and limited. \cite{nguyen2021temporal} introduced a latent space model to capture the complex distribution of the input data. The uncertainty of the time-series input is learned through a probabilistic latent space.

\subsection{Decision-making under uncertainty}

The most prevalent methods for addressing decision-making in uncertain conditions is \gls{abk:so} \citep{birge2011introduction}. \gls{abk:so} aims to identify robust and near-optimal solutions in scenarios with uncertain objectives or constraints. An alternative strategy is \gls{abk:ro}. Unlike \gls{abk:so}, which relies on a known or partially known uncertainty distribution, \gls{abk:ro} requires only a range of possible values for the uncertain parameters \citep{ben2009robust}. \gls{abk:ro} has been extensively employed in solving optimization challenges in large networks where uncertainty is a constant factor, such as in network facility location \citep{baron2011facility,gulpinar2013robust} and inventory management \citep{bienstock2008computing,see2010robust}.

Recently, \gls{abk:dro} \citep{Rahimian_2022} has gained traction as a method for decision-making under uncertainty. Distinct from the \gls{abk:ro} framework, which heavily relies on the accurate definition of the uncertainty set, \gls{abk:dro} incorporates probability information into the optimization process and provides more flexibility in defining uncertainty. Specifically, \gls{abk:dro} aims to establish a collection of potential distributions. The set of potential distributions is usually called an \emph{ambiguity set}. This set is typically formulated using data-driven methods, often characterized by moment or cross-moment constraints. \gls{abk:dro} has found its applications in a variety of dynamic system-related optimization problems, including portfolio selection \citep{delage2010distributionally} and scheduling medical appointments \citep{bertsimas2019adaptive}.

\subsection{\gls{abk:amod} vehicle rebalancing}

Since their appearance, on-demand systems have been extensively investigated by previous studies \cite{wang2019ridesourcing}. A critical element of the operations in these systems is the rebalancing of idle vehicles, as rebalancing is an essential strategy to ensure efficient allocation of drivers to meet customer demand when a mismatch exists between the demand and supply. Many previous studies have concentrated on using operations research techniques for rebalancing vehicles in on-demand services\citep{spieser2016shared, 8593743, pavone2012robotic,zhang2014control, 10.1145/3055004.3055024, 8317908, GUO2021161, 9744407}. By employing optimization models, a wide range of objectives and assumptions can be incorporated through the adjustment of objective functions and the integration of various constraints.


In the optimization-based rebalancing models, demand prediction is a crucial input for deriving optimal operations. While many existing studies assume a known and deterministic demand as discussed in \cref{demand}, incorporating uncertainty in demand is more realistic and practical in real-world settings. Demand prediction and rebalancing strategies have long been central to vehicle rebalancing research, but there are limited studies exploring the synergy of the two approaches. Considering that demand information is crucial for designing rebalancing schemes, and demand prediction primarily serves to facilitate system management, integrating these two aspects is particularly relevant for real-world applications. 

Among the limited research in this area, \cite{ALKANJ20201088} focused on a series of decisions in electric vehicle systems, including car dispatching (recharging, repositioning, parking), pricing, and fleet sizing. Their approach utilizes a value function to assess decision quality in the face of demand uncertainty, where the uncertainty is not pre-forecasted but captured by the value function during system state transitions. \cite{RAMEZANI2018203} introduced a predictive control method for taxi dispatching, using a Macroscopic Fundamental Diagram (MFD) graph to update system states and form the basis for prediction and optimization, incorporating uncertainty by adding Gaussian noise to point predictions. Meanwhile, \cite{10018486} proposed a distortion function and a dual distortion function to quantify the uncertainty of cumulative rewards, linking this uncertainty to vehicle actions rather than demand distribution.

Another effective strategy to manage demand uncertainty in on-demand systems is the "prediction-then-optimize" approach. This method can be described by the following process: the uncertainty is predicted and quantified using historical data, and the prediction is used as input to the optimization model. Following this paradigm, \cite{RN19} adopted a hypothesis testing-based method to develop uncertainty sets from demand data. These sets are then applied to address the vehicle rebalancing problem through \gls{abk:ro}. Expanding on this idea, the same researcher further evolved the technique to create ambiguity sets for use in \gls{abk:dro}~\cite{RN20}. Additionally, \cite{9341481} introduced a framework tailored for electric vehicle system rebalancing under both demand and supply uncertainties. This framework leveraged uncertainty sets derived from the AutoRegressive Integrated Moving Average (ARIMA) model \citep{ZHANG2003159}. This model was subsequently incorporated into a \gls{abk:dro} model. The proposed model considered not only vehicle rebalancing but also the utilization of charging stations. However, it did not account for vehicle-request matching.

A recent study \cite{9744407} presented an innovative predictive prescription method. This approach leverages historical demand data, applying unsupervised machine learning models to assign probabilities to each historical demand scenario. These probabilities are informed by auxiliary data such as weather conditions, temperature, and time of day. The method then employs a \gls{abk:so} model, which integrates the historical data and the derived probabilities to formulate rebalancing decisions. However, this approach has its limitations. It requires a limitation on the number of demand scenarios for the \gls{abk:so} to be solvable, since each additional demand scenario introduces new constraints and variables into the model. As a result, an extensive demand set can lead to a model that becomes too complex to manage effectively. On a different note, \cite{GUO2023104244} proposed an alternative approach involving a data-driven robust optimization model. This model is unique in that it constructs its uncertainty set through a series of point prediction models, each subject to a defined structural risk. A key limitation of this model is that the structural risk of the prediction model must be a convex function, which limits the range of applicable models.

To bridge the above research gap, in this paper, we introduce a novel framework that synergizes \gls{abk:dnn}s with \gls{abk:ro} to address the \gls{abk:amod} rebalancing challenge. The goal is to harness the power of \gls{abk:dnn} in uncertainty quantification for uncertain parameters defined in \gls{abk:ro}, thereby enhancing the performance of the \gls{abk:ro} model by better constructing the uncertainty set. \gls{abk:dnn}s are renowned for their ability to effectively discern complex relationships between observed and target datasets, offering a reliable estimation of the uncertainty set. Consequently, this integrated optimization framework is poised to deliver superior performance in practical, real-world scenarios.
\section{Methodology}\label{method}

In this section, we propose a general \gls{abk:duro} framework. This framework is adaptable for any decision-making problems with uncertainty. Subsequently, we present the details of the proposed upstream demand prediction module and the downstream optimization module for vehicle rebalancing.

\subsection{Framework Overview}
\begin{figure}[tb]
    \centering
    \includegraphics[width=0.5\textwidth]{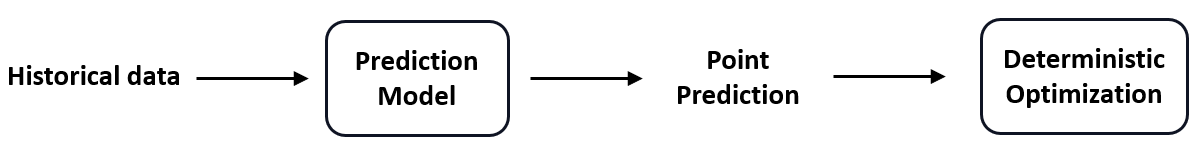}
    \caption{Common workflow for solving the rebalancing problem. 
    The prediction provides point demand prediction based on historical data and the predicted values are used as input to a downstream deterministic optimization model.}
    \label{fig:f_point}
\end{figure}

In the prevailing literature, the optimization problem in dynamic systems is conventionally approached through the framework shown in \cref{fig:f_point}. Various studies are concentrating on refining point prediction accuracy of the decisive variables, such as the demand in vehicle rebalancing case, and others on developing optimization models and enhancing algorithms. Although both prediction and optimization methodologies have individually demonstrated success in delivering accurate point prediction and improving system performance, a composition of these two components does not necessarily ensure a satisfactory result, primarily due to the pervasive uncertainty inherent in real-world scenarios.

\begin{figure}[tb]
    \centering
    \includegraphics[width=0.5\textwidth]{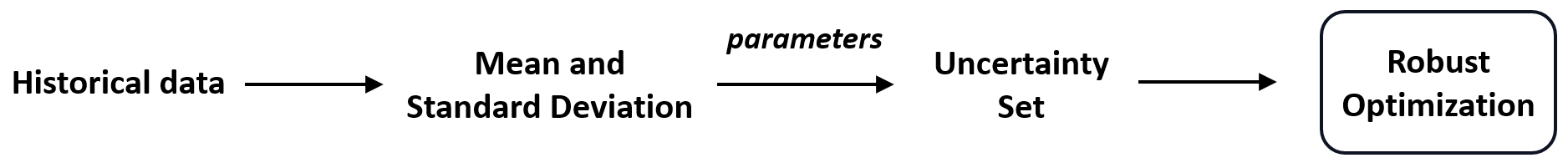}
    \caption{Optimization workflow proposed by~\cite{9744407}. The mean and standard deviation is retrieved from historical data and utilized to compute the uncertainty set for the robust optimization model.}
    \label{fig:f_ro}
\end{figure}

In the work presented by \cite{9744407}, the authors proposed a momentum-based robust optimization model to address an optimization problem in the presence of uncertainty. The framework is outlined in \cref{fig:f_ro}. The uncertainty set contains all possible values of the random variables. This framework is defined based on the mean and standard deviation calculated from the historical data and additional parameters. The uncertainty set is subsequently input to the robust optimization model. The robust optimization model solves the optimization problem against all possible demand patterns.

A notable drawback of the outlined framework is the manual selection of parameters used to define the uncertainty set. These parameters lack a data-driven foundation and are based on the operators' subjective interpretations, necessitating a substantial amount of parameter tuning. Selecting these parameters lacks standardization, introducing inefficiency into the framework. 

\begin{figure}[tb]
    \centering
    \includegraphics[width=0.4\textwidth]{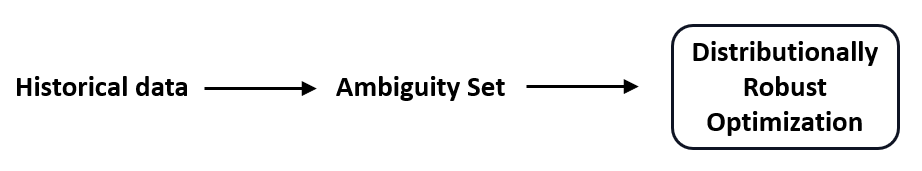}
    \caption{\gls{abk:dro} workflow. The historical data is used to build a set of potential distributions of uncertain variables, called the ambiguity set, and subsequently used in the distributional robust optimization model.}
    \label{fig:f_dro}
\end{figure}

Another commonly used model for handling uncertainty is \gls{abk:dro} \citep{10.1145/3055004.3055024}. The workflow of \gls{abk:dro} is presented in \cref{fig:f_dro}. Compared to \gls{abk:ro} that quantifies uncertainty through an uncertainty set containing possible variable values, \gls{abk:dro} defines the potential distributions of the random variable within an ambiguity set. The ambiguity set is derived from data and does not depend on tuned parameters. However, since the ambiguity set usually involves a large number of distributions to test, the \gls{abk:dro} model embeds a higher complexity than the \gls{abk:ro} model and is less computationally efficient. This limitation underscores the need for an automated and data-driven approach to enhance the robustness and efficiency of the uncertainty quantification process in vehicle rebalancing problems, which is essential to solve the problem in real time.

\begin{figure}[tb]
    \centering
    \includegraphics[width=0.4\textwidth]{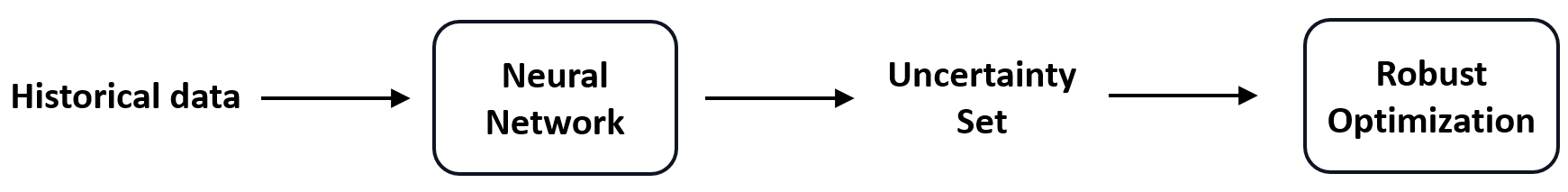}
    \caption{DURO workflow. A neural network is leveraged to predict the uncertainty set used in the robust optimization model.}
    \label{fig:general}
\end{figure}

To mitigate the limitations of existing methods, we propose an innovative neural network-based robust optimization framework, \gls{abk:duro}, as shown in \cref{fig:general}. This framework is generalizable, applicable to various optimization problems arising in systems marked by uncertainty. It operates as an end-to-end system, seamlessly progressing from data to decision. In \gls{abk:duro}, we employ a deep neural network to predict the uncertainty set using historical data, as opposed to manually crafting it as seen in \cref{fig:f_ro}. Compared to \gls{abk:dro}, this approach ensures enhanced computational efficiency and thus is more applicable for real-world scenarios.

\begin{figure}[tb]
    \centering
    \includegraphics[width=0.5\textwidth]{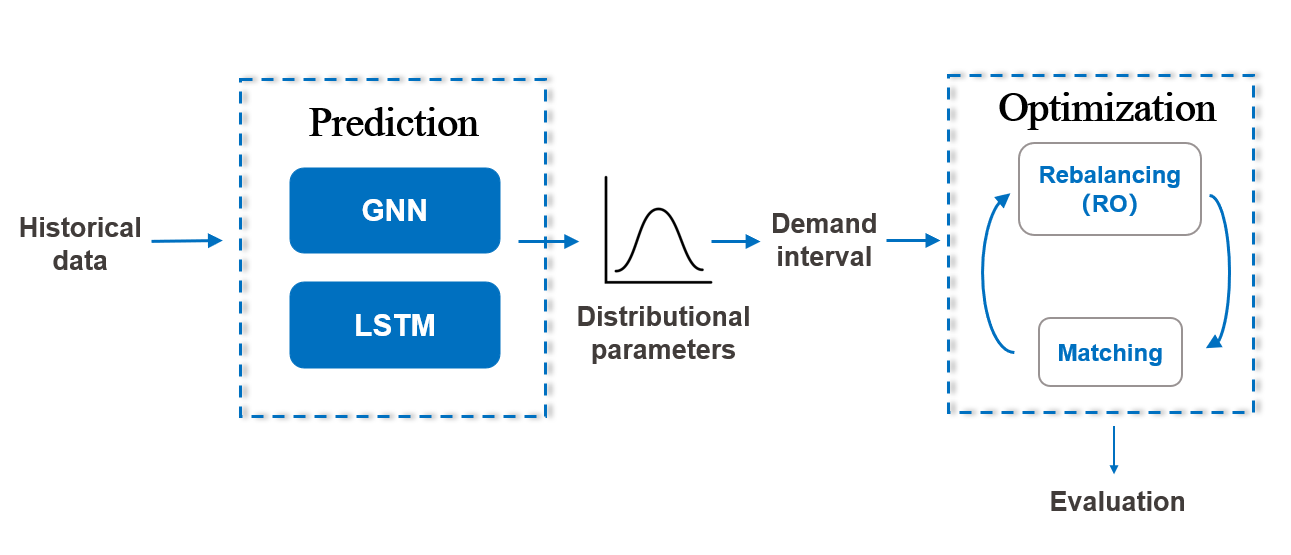}
    \caption{Framework architecture. There are two modules: prediction and optimization. The prediction module is for demand uncertainty quantification and the robust optimization model is incorporated in the optimization module to solve the rebalancing problem.}
    \label{framwork}
\end{figure}

In the rest of this paper, we demonstrate the efficacy of \gls{abk:duro} by applying it to the vehicle rebalancing problem in \gls{abk:amod} systems and comparing it to other aforementioned methods. \cref{framwork} presents an overview of the vehicle rebalancing framework architecture based on \gls{abk:duro}. The prediction module is a neural network primarily composed of a graph neural layer and an \gls{abk:lstm} layer. The graph neural layer encodes the spatial information in the demand dataset, and the \gls{abk:lstm} captures the temporal information. The output of this network comprises the parameters relating to the probabilistic assumptions governing the demand distribution. Subsequently, these parameters are utilized to compute the demand interval, which is integrated into the following optimization module.

The optimization module functions as a simulator that encompasses both rebalancing and matching. Rebalancing and matching are executed iteratively until a predefined number of iterations is reached. The evaluation metrics are then calculated based on the simulation's outcomes. This framework aims to showcase the holistic performance of \gls{abk:duro} in comparison to other existing methodologies.

In the following subsections, we formally present the details of its two modules: the uncertainty quantification deep neural network and the robust vehicle rebalancing model. 
 
\subsection{Demand uncertainty quantification}\label{uncertainty}
\begin{figure}[tb]
    \centering
    \includegraphics[width=0.5\textwidth]{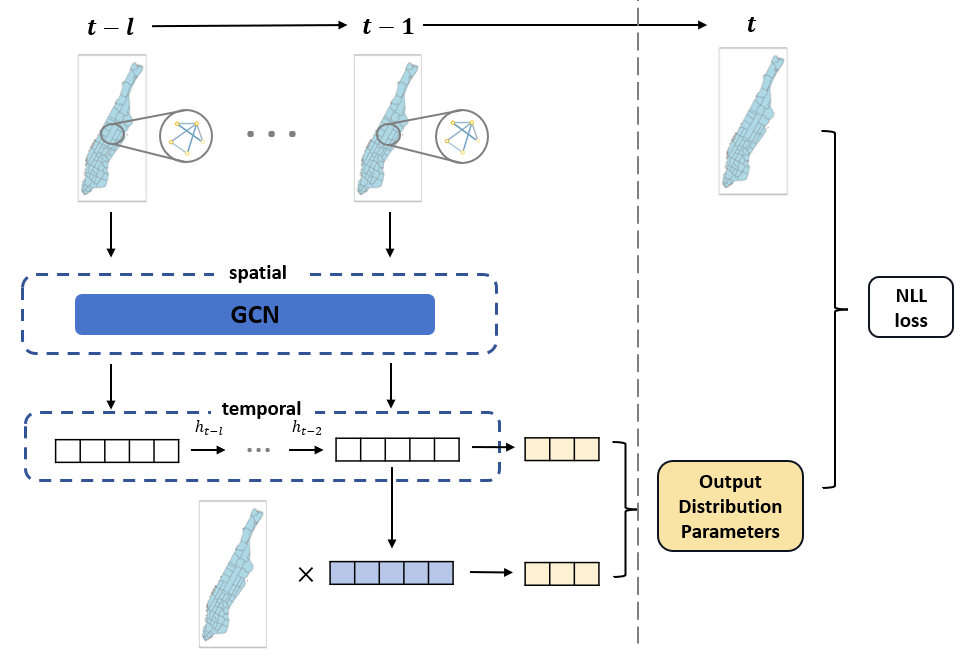}
    \caption{Proposed model architecture. The model encodes past demand information by spatial and temporal layers. The processed information is added with auxiliary information to estimate the distributional parameters by negative log-likelihood (NLL) loss.}
    \label{nn}
\end{figure}

This section provides a detailed description of the proposed uncertainty prediction module. The prediction module deals with the demand forecasting problem under uncertainty by a probabilistic graph convolutional neural network characterized by the following formula:
\[y \sim G(\theta) = G(F(\boldsymbol{X},w)),\]
where $G(\theta)$ is the probabilistic assumption about $y$, which is the demand, and $F(\boldsymbol{X},w)$ parametrizes the probability distribution. In our case, $F(\boldsymbol{X},w)$ is represented by a neural network. \cref{nn} provides an overview of the architecture of the neural network.
The neural network mainly comprises graph convolutional layers and \gls{abk:lstm} layers. \gls{abk:gcn} has been widely applied and has been shown to be effective in predicting transportation demand \citep{JIANG2022117921}. It is commonly used to encode the spatial information in a network, updating each node's feature while considering its related neighbor nodes and arcs. For each \gls{abk:gcn} layer $l$, we follow the following propagation rule:
\[H^{l+1} = \sigma(\widetilde{D}^{-\frac{1}{2}}\widetilde{A}\widetilde{D}^{-\frac{1}{2}}H^{l}W^{l}),\]
where $\widetilde{A}$ is the adjacency matrix, $\widetilde{D}$ is a diagonal degree matrix and $W^{l}$ is the layer-specific trainable weight matrix. $\sigma(\cdot)$ denotes the activation function and we use $ReLU$ here.

The encoded spatial features are provided as input to the \gls{abk:lstm} layers. \gls{abk:lstm} is designed for dealing with time-series data with a mechanism to treat long and short-term information differently. In each \gls{abk:lstm} cell, the temporal correlation is learned and captured by the forget, input, and output gates. The encoded long and short-term memory, which are denoted as cell state and hidden state, maintains the overall past information and is passed to the next timestamp. 

The output of the \gls{abk:gcn}-\gls{abk:lstm} layers is subsequently passed to two components: a linear layer for spatial-temporal prediction of the recent demand and another linear layer to connect the influence of the historical average demand. The two components are summed up to generate the final prediction. The final prediction is the parameters of the distributional assumptions made on the demand distribution. With the parameters of the corresponding parameters as the output from the network, the probability of the true ridership can be sampled, and this probability is used to train the complete network by minimizing negative log-likelihood (NLL):
\[NLL = - \sum_{i} G(y_{i}|X_{i};w),\]
\subsection{Matching-Intergrated Vehicle Rebalancing Model} \label{ro model}
We propose a robust optimization model that optimizes the rebalancing operations while accounting for matching. The rebalancing problem is solved following the paradigm of \gls{abk:mpc} to ensure non-myopic decisions. Before getting into the details of the full model, we first clarify the problem setting and some notations that we will use in the following explanation. 

\begin{figure}[tb]
    \centering
    \includegraphics[width=0.5\textwidth]{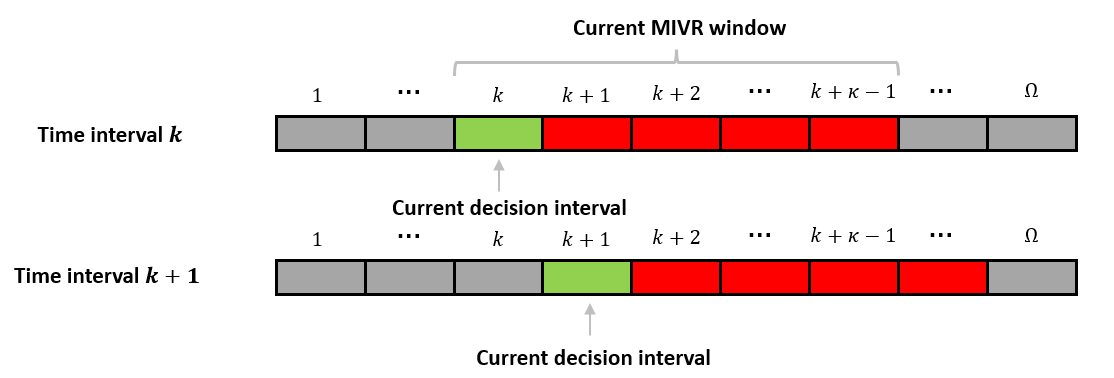}
    \caption{Example of rolling optimization. The green interval indicates the current decision time interval and the red ones are the future intervals that are also considered in the current optimization horizon.}
    \label{mivr}
\end{figure}

Suppose that we have identified an operation period of interest, we first divide it into $\Omega$ identical time intervals index by $k = 1,2, ..., \Omega$, with an interval length $\Delta$. The MIVR model is solved on a rolling basis, which means the decision variables are determined repeatedly at the beginning of each time interval. Note that even though we only implement the decisions of the current interval, when solving the model for a time interval $k$, we still incorporate the optimization problem for the $\kappa$ future time intervals, as the operation for the current time interval is highly dependent on the future situation. This process is displayed in Figure \ref{mivr}. The space is discretized into $n$ regions index by $i = 1, ..., n$. Combined with the time index, the demand of region $i$ at time $k$ is denoted by $r_i^k$. We have sets $N = \{1, .., n\}$ for all the region and $K = \{1, ..., \kappa\}$ for all time intervals. 
For each time interval, we do both rebalancing and passenger-vehicle matching. For vehicle rebalancing, we assume it takes place at the beginning of each time interval and the matching happens at the end of it. Next, we will present the core constraints for each phase.

In the rebalancing phase, we denote the number of idle vehicles from region $i$ to $j$ at time $k$ by $x_{ij}^{k}$. Let $S_i^k$ indicate the number of available vehicles in region $i$ at the end of time $k$. $d_{ij}^{k}$, $w_{ij}^{k}$ denote the travel distance and time from region $i$ to $j$ at time $k$ respectively. The values of these two sets of parameters are approximated by the distance and travel time between the centroids of the two regions. We define a variable $a_{ij}^k \in \{0, 1\}$ denoting whether an idle vehicle can be rebalanced from region $i$ to $j$ at time $k$ where $a_{ij}^k = 0$ indicates the rebalancing is feasible. $a_{ij}^k = 0$ if $w_{ij}^k \leq \Delta $. This feasibility constraint is given by:
\begin{equation}
    a_{ij}^{k} \cdot x_{ij}^{k} = 0 \quad \forall i, j \in N, \forall k \in K \label{eq:1},
\end{equation}
When solving the rebalancing model, the actual quantity and location of each passenger and vehicle are not yet available, so the matching component is solved based on the estimated demand. Denote the decision variables for the matching phase as $y_{ij}^k$, which indicates the number of customers in region $i$ matched with vehicles in region $j$ at time $k$. Let $b_{ij}^k \in \{0, 1\}$ denote the feasibility of the matching, where $b_{ij}^k = 0$ indicates feasible. $b_{ij}^k = 0$ only if the pickup time is less than an imposed maximum value. Let $\bar w$ denote the maximum allowed pickup time, the feasibility constraint is:
\begin{equation}
    b_{ij}^{k} \cdot y_{ij}^{k} = 0 \quad \forall i, j \in N, \forall k \in K \label{eq:2}
,\end{equation}
It is worth mentioning that the matching component only serves as the auxiliary component in the model, which focuses on optimizing the rebalancing decisions. In a real-world scenario, when the vehicle is rebalanced and requests are collected, the matching problem can then be solved by a separate dedicated matching model given the realized demand. This is also how we run the simulation in the experiments.

The number of matched trips is upper bounded by the available vehicles and the demand. These constraints are:
\begin{subequations}
\begin{align}
    \sum_{j=1}^{n} y_{ji}^{k} \leq S_i^k \quad \forall i \in N, \forall k \in K \label{eq:3a},\\
    \sum_{j=1}^{n} y_{ji}^{k} \leq r_i^k \quad \forall i \in N, \forall k \in K \label{eq:3b},\\
    T_i^k = r_i^k - \sum_{j=1}^{n} y_{ij}^k \quad \forall i \in N, \forall k \in K \label{eq:3c}.
\end{align}
\end{subequations}
Constraints (\ref{eq:3a}), (\ref{eq:3b}) are the upper bounds constraints, and $T_i^k \in \mathbb{R_+}$ in constrain (\ref{eq:3c}) represents the number of unsatisfied demand.

The rebalancing and matching phases are connected by the vehicle transition process in each region. To establish this relationship, we define more variables related to the number of available and occupied vehicles:
\begin{itemize}
    \item $V_i^k$: Number of available vehicles in region $i$ at the beginning of time interval $k$.
    \item $O_i^k$: Number of occupied vehicles in region $i$ at the beginning of time interval $k$.
    \item $P^k(P_{ij}^k)$: The probability that an occupied vehicle located in region $i$ at time $k$ will be in region $j$ and stay occupied at time $k+1$.
    \item $Q^k(Q_{ij}^k)$: The probability that an occupied vehicle located in region $i$ at time $k$ will be in region $j$ and become vacant at time $k+1$.
\end{itemize}

$P^k$, $Q^k$ are regional transition matrices describing the movement of occupied vehicles. We approximate them with static matrices estimated from historical data. Interested readers can refer to \cite{GUO2021161} for more details and discussion on this method.

With the above variables, the rebalancing and matching phases can be connected by the following constraints:
\begin{subequations}
\begin{align}
    & \sum_{j=1}^{n} x_{ji}^{k} \leq V_i^k \quad \forall i \in N, \forall k \in K ,\label{eq:4a}\\
    & S_i^k = V_i^k + \sum_{j=1}^{n} x_{ji}^k - \sum_{j=1}^{n} x_{ij}^k \quad \forall i \in N, \forall k \in K, \label{eq:4b}\\
    & V_i^{k+1} = S_i^k - \sum_{j=1}^{n} y_{ji}^k + \sum_{j=1}^{n} Q_{ji}^{k}O_{j}^k \quad \forall i \in N, \forall k \in K \backslash \{\kappa\}, \label{eq:4c}\\
    & O_i^{k+1} = \sum_{j=1}^{n} y_{ji}^k + \sum_{j=1}^{n} P_{ji}^{k}O_{j}^k \quad \forall i \in N, \forall k \in K \backslash \{\kappa\}, \label{eq:4d}
\end{align}
\end{subequations}    
Constraint (\ref{eq:4a}) ensures that the total number of vehicles rebalanced from region $i$ to other regions at time $k$ does not exceed the available vehicles at the beginning of the time interval. Constraint (\ref{eq:4b}) states that the gap between the number of available vehicles at the beginning and end of a time interval in one region is equal to the difference between the vehicles rebalanced to and from the region. Constraint (\ref{eq:4c}) shows that the available vehicles at the beginning of the time interval $k+1$ is related to the available number of vehicles at the end of time $k$, the matched trips at time $k$, and the number of occupied vehicles that become available at the beginning of time $k+1$. Constraint (\ref{eq:4d}) states a similar relationship for the number of occupied vehicles. 

The objective of the model is to minimize the number of unsatisfied trips and the total vehicle traveling distance, represented in the following formula:
\begin{dmath}
    \min \quad Z = \sum_{k=1}^{\kappa}\sum_{i=1}^{n}\sum_{j=1}^{n} x_{ij}^{k}d_{ij}^{k} + \beta \cdot \sum_{k=1}^{\kappa}\sum_{i=1}^{n}\sum_{j=1}^{n} y_{ij}^{k}d_{ij}^{k} + \gamma \cdot \sum_{k=1}^{\kappa}\sum_{i=1}^{n} T_{i}^{k}, 
\end{dmath}
where $\beta$ is a parameter that defines the relative weighting of rebalancing distance and pick distance, and $\gamma$ indicates the penalty induced by each unsatisfied trip. 

In the model, we consider the uncertainty of the demand variable $r_i^k$. In real-world applications, the future demand is unknown during the planning phase and is usually predicted by some model. As the prediction can hardly be exactly the same as the true value, we should consider the prediction's uncertainty in the optimization model. For MIVR, we define the uncertainty set for demand $r_i^k$ based on the predicted  $\mu_i^k$ and interval $(LB_i^k, UB_i^k)$. The box uncertainty can be imposed by the following uncertainty constraints:
\begin{subequations}
\begin{align}
    & LB_i^k \leq r_i^k \leq UB_i^k \quad \forall i \in N, \forall k \in K, \label{eq:6a}\\
    & |\sum_{i=1}^{n} (r_i^k - \mu_i^k) | \leq \Gamma \quad \forall k \in K, \label{eq:6b}
\end{align}
\end{subequations} 
Constraint (\ref{eq:6a}) is a box uncertainty set that imposes the upper and lower bounds from the estimated region interval at each time interval $k$. (\ref{eq:6b}) is a polyhedral uncertainty set that limits the total offset in the sum of the demand during a time interval across all regions. To avoid equality constraints with uncertain parameters facility the computation, the robust optimization model can be reframed into the following format:
\begin{dmath}
    \min \quad Z = \sum_{k=1}^{\kappa}\sum_{i=1}^{n}\sum_{j=1}^{n} x_{ij}^{k}d_{ij}^{k} + \beta \cdot \sum_{k=1}^{\kappa}\sum_{i=1}^{n}\sum_{j=1}^{n} y_{ij}^{k}d_{ij}^{k} + \gamma \cdot \sum_{k=1}^{\kappa}\sum_{i=1}^{n} (r_i^k - \sum_{j=1}^{n} y_{ij}^k),     
\end{dmath}
with constrains (\ref{eq:1}), (\ref{eq:2}), (\ref{eq:3a}), (\ref{eq:3b}), (\ref{eq:4a}) - (\ref{eq:4d}).

\section{Experiments}\label{experiment}
In this section, we present prediction and simulation results of the proposed framework based on real-world on-demand data from New York City. The goal of the experiments is to address the following questions: (1) What is the most proper probabilistic assumption to make to get the best uncertainty set? (2) How does \gls{abk:duro} perform compared to other optimization approaches? To answer the second question, we conduct simulations with different optimization methods. The tested methods encompass a range of demand prediction models and optimization frameworks.
\subsection{Data Description}
In the experiments, we use an open-source on-demand dataset for the island of Manhattan in New York City. It includes trip information from businesses that dispatch or plan to dispatch more than 10,000 trips in New York City per day, including Uber, Lyft, Via, and Juno. In the dataset, the Manhattan area is discretized into 63 taxi zones, and each zone is assigned a unique zone ID. Trip information includes the pickup and dropoff time, along with corresponding zone IDs. Our experiments employ the trip data on all weekdays in June 2019, which is 20 days in total.

\begin{figure}[h]
    \centering
    \begin{subfigure}[t]{0.2\textwidth}
        \centering
        \includegraphics[width=0.7\textwidth]{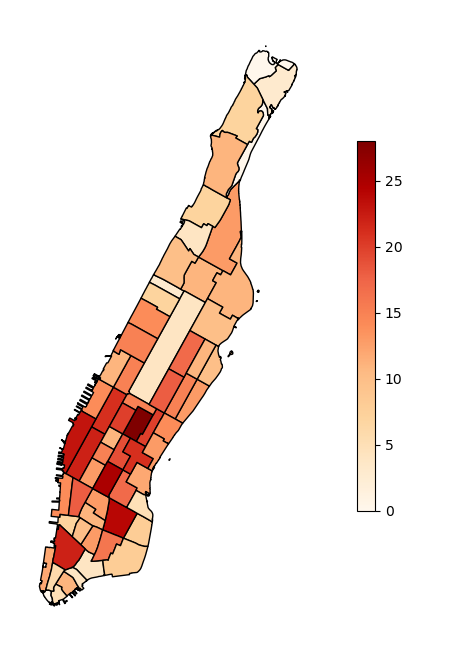}
        \caption{Average demand}
    \end{subfigure}%
    ~ 
    \begin{subfigure}[t]{0.25\textwidth}
        \centering
        \includegraphics[width=1.0\textwidth]{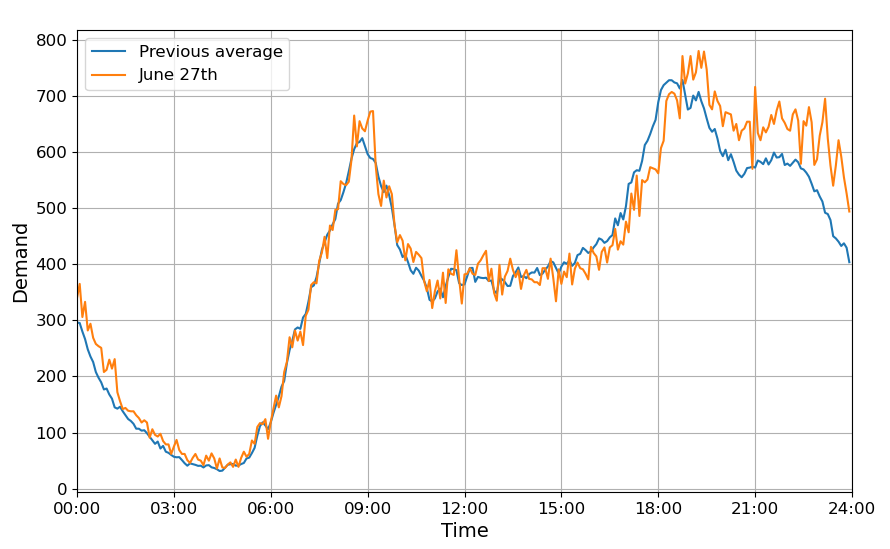}
        \caption{Total demand fluctuation}
    \end{subfigure}
    \caption{Demand pattern of NYC Manhattan. The demand is unevenly distributed across different regions as shown in (a) and fluctuates during the day as in (b). (b) also shows that the demand pattern is similar on the same day of the week in June as the historical average is closely aligned with that on the 27th.}
    \label{fig:demand}
\end{figure}

During this period, there are on average around 20,000 trips per day. For each origin-destination pair, We aggregate all the trips in 5-minute intervals. Figure \ref{fig:demand} (a) displays the average demand on June 26th, 2019. This figure reveals an imbalance of the demand across different zones, with a notable demand concentration in the southern region. This disparity also underscores the necessity of rebalancing to serve the non-uniformly distributed trip demands. In Figure \ref{fig:demand} (b), the total demand fluctuation on the same day is juxtaposed with the average demand from previous days. The morning peak hour is observed around 9:00 and the night peak around 19:00. The resemblance in the shapes of these two lines suggests a recurring demand fluctuation pattern within this region.
\subsection{Evaluation and Simulator}
To assess the overall performance of the complete DURO framework, we first select the most suitable probabilistic assumption based on its performance. As discussed in \cref{uncertainty}, the network is trained using NLL. To comprehensively evaluate the quality of different distributional assumptions, additional metrics are introduced alongside NLL. The utilized metrics can be classified into three categories: composite (NLL), point, and uncertainty. 

For point prediction accuracy, Mean Absolute Error (MAE) and Mean Absolute Percentage Error (MAPE) are employed. Regarding uncertainty, two key indicators include Prediction Interval Coverage Probability (PICP) and Mean Prediction Interval Width (MPIW). PICP indicates the percentage of the actual data points covered by the prediction interval, while MPIW quantifies the average width of the prediction intervals. MPIW is defined as follows\citep{5672788}:
\[MPIW = \frac{1}{n} \sum_i^n(U_i - L_i),\],
where $U_i$ and $L_i$ represent the upper and lower bounds of the prediction interval, and $n$ is the number of predicted intervals. The intended coverage percentage we use for MPIW and PICP to compare the different distributions is 95\%.

\begin{figure}[tb]
    \centering
    \includegraphics[width=0.5\textwidth]{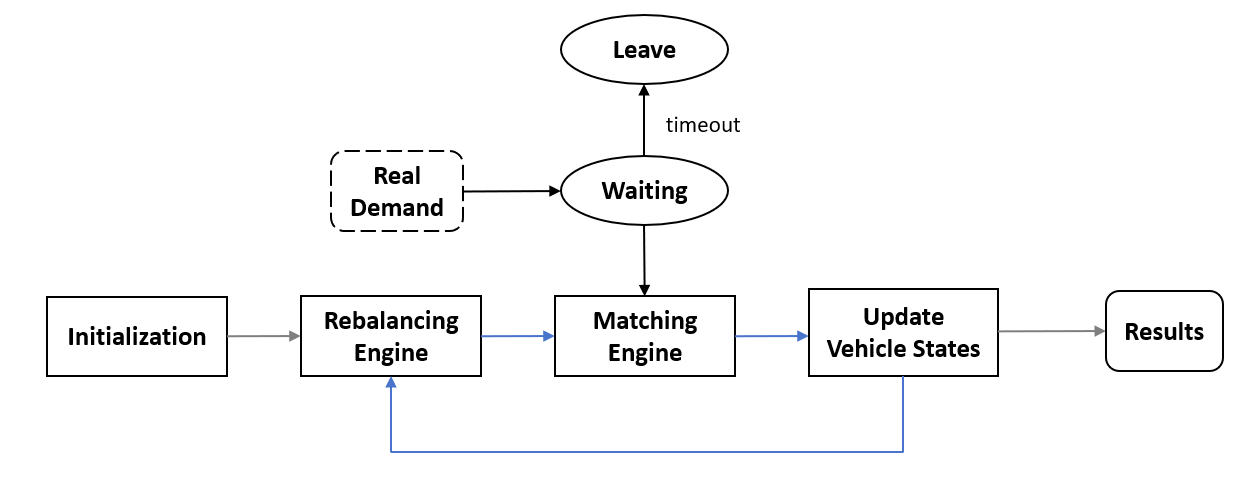}
    \caption{Simulation workflow in the experiments. At each time step, the rebalancing engine computes the rebalancing decisions and the matching engine matches available vehicles with passengers. The vehicles move towards the rebalancing destination or pickup points of the passengers, and then system status is updated for the next decision window.}
    \label{simulator}
\end{figure}

The performance of the complete framework is evaluated by running simulations. The workflow of the simulator is outlined in \cref{simulator}. The initialization step includes random initialization of the vehicle locations and loading the road networks for the Manhattan area with the shortest path distance, the predecessor matrices, travel time matrices between taxi zones, and regional transition matrices. During each simulation step, the rebalancing decisions are formulated and implemented at the beginning of the step, and the matching between passengers and vehicles is executed at the end. The matching engine works based on the passengers generated from real demand data. Each passenger has a maximum waiting time of 300 seconds and will leave the system immediately after reaching this maximum.

During the simulation, we recorded the served passenger's waiting time and traveling time. The waiting time is defined as the duration a passenger spends waiting from the moment they send a request until they are picked up, and the traveling time represents the time it takes for the vehicle to travel from the origin to the destination. Upon completion, an overall passenger leaving rate was computed. These metrics serve as crucial indicators of the performance of different rebalancing engines. The evaluated rebalancing engines will be introduced in the subsequent benchmark subsection. 
\subsection{Benchmarks}
The performance evaluation of \gls{abk:duro} involved a comparative analysis with various optimization models, wherein the rebalancing engine was replaced with the corresponding model. In the following sections, these optimization models are referred to as different rebalancing engines. \gls{abk:duro} is benchmarked against the following optimization engines:
\begin{enumerate}
    \item \textbf{\emph{\gls{abk:dohv}}}: In this rebalancing engine, the demand is predicted based on the average demand before June 27. The prediction is assumed to be accurate and utilized in the deterministic version of the optimization model outlined in \cref{ro model}, excluding the demand uncertainty constraints.
    \item \textbf{\emph{\gls{abk:donn}}}: In contrast to the historical average demand approach, \gls{abk:donn} predicts the demand by a graph-\gls{abk:lstm}-based neural network. This neural network is composed of two graph convolutional layers followed by a \gls{abk:lstm} layer. The predicted demand is subsequently fed into the same optimization model as described in engine 1.
    \item \textbf{\emph{\gls{abk:ro}}}: The \gls{abk:ro} engine leverages the same robust optimization model as in the \gls{abk:duro}, which is presented in \cref{ro model}. The difference between this engine and \gls{abk:duro} is that, for constraint (\ref{eq:6a}), \gls{abk:duro} uses a neural network to predict the lower and upper bound, while in \gls{abk:ro} engine, this constraint is presented in the following format:
    \[-\rho \leq \frac{r_i^k - \mu_i^k}{\sigma_i^k} \leq \rho,\]
    where $\mu_i^k$, $\sigma_i^k$ is the mean and standard deviation of the demand data before June 27, and $\rho \in \mathbb{R^+}$. $\rho$ indicates how far the real demand can be away from the historical average in proportion to the historical standard deviation.
    \item \textbf{\emph{\gls{abk:dro}}}: This engine applies the distributional version of the model presented in \cref{ro model}. Instead of specifying uncertainty set by Constraints (\ref{eq:6a}) and (\ref{eq:6b}), the uncertainty is indicated by ambiguity set. The ambiguity set is captured by the cross-moment ambiguity deriving from historical data \citep{delage2010distributionally, bertsimas2019adaptive, wiesemann2014distributionally,chen2020robust}:
    \begin{gather*} 
    \boldsymbol{\widetilde d} \sim \mathbb{R}\\ 
    \mathbb{E_P}[\boldsymbol{\widetilde d}] = \boldsymbol{\mu} \\
    \mathbb{E_P}[(\widetilde d_j - \mu_j)^2] \leq \sigma_j^2 \\
    \mathbb{E_P}[\boldsymbol{e}^T(\boldsymbol{\widetilde d} - \boldsymbol{\mu})^2] \leq \boldsymbol{e}^T \Sigma \boldsymbol{e} \\
    \mathbb{P}[\boldsymbol{\widetilde d} \in [\boldsymbol{\underline d}, \boldsymbol{\bar d}]] = 1,
    \end{gather*}
\end{enumerate}
for all $j \in N$, where $\boldsymbol{\widetilde d}$ is the uncertain variables, $\boldsymbol{\mu}$ is the historical mean vectors, $\boldsymbol{\sigma}$ is the historical standard deviations vector, and $\Sigma$ is the covariance matrix.

Table \ref{parameters} displays the model parameters shared across all the optimization engines in the simulation.
\begin{table*}[h]
\centering
\begin{tabular}{|l|l|l|}
\hline
Parameter & Description & Value\\ \hline
 $\beta$ & Weigh parameter for pickup distance & 1\\ \hline
 $\gamma$ & Penalty for unsatisfied requests & 100 \\ \hline
 $\Omega$ & Simulated time intervals & 24 \\ \hline
 $\Delta$ & Decision time interval length for vehicle rebalancing problem & 300 s \\ \hline
 $\delta$ & Decision time interval length for matching problem & 30 s \\ \hline
 $\overline{\omega}$ &Maximum pickup time & 30 s \\ \hline
 $\widetilde{\omega}$ & Maximum passenger waiting time & 300 s \\ \hline
 $n$ & Number of regions & 63 \\ \hline
 $\kappa$ & Number of look-ahead time intervals when solving MIVR model & 6 \\ \hline
 $m$ & Number of historical data points & 18 \\ \hline
 $N_v$ & Number of vehicles & 2000 \\ \hline
\end{tabular}
\caption{Parameter values used in simulation}
\label{parameters}
\end{table*}


\subsection{Results}

In this section, we present the results of deep uncertainty quantification and simulations. Based on the simulation results, we provide a comparative analysis of \gls{abk:duro} against other optimization engines. To train the demand uncertainty quantification network, we partitioned the dataset into training, validation, and test sets along the temporal dimension. The training set contains data from workdays between June 1, 2019, and June 24, 2019, the validation set includes data from June 25 and 26, and the test set comprises data from June 27. We experimented with various probabilistic assumptions and determined that the Poisson distribution yielded the best accuracy. Consequently, the Poisson distribution is employed in the prediction module of DURO for subsequent performance testing.

The different engines were evaluated by running simulations based on the real data from the morning peak between 7:00 -9 :00 on June 27, spanning a total number of 24 time intervals. When utilizing the rebalancing engine, it is assumed that all evaluated engines have access to the data before June 27. The actual demand on June 27 is treated as future demand and thus remains unknown when running the optimization models. It is only used for generating passengers during the matching phase in the simulator. The results of the different engines are summarized in \cref{prediction}. For engines with various parameter combination cases, we show the result of the best-performing one.
\begin{table}[h]
\centering
\begin{tabular}{|c|c|c|c|}
\hline
\textbf{Engine} & \textbf{Average Waiting (s)} & \textbf{Leaving Rate (\%)} \\ \hline
DOHV   &  127.3                            &    5.1                       \\ \hline
DONN   &  133.3                            &     4.9                      \\ \hline
RO   &  123.8                            &      \textbf{4.8}                     \\ \hline
DRO*     &   133.5                           &    5.3                       \\ \hline
DURO   &   \textbf{123.5}                          &        4.9                   \\ \hline
\end{tabular}
\caption{Performance of different engines. $\boldsymbol{DRO^*}$: The rebalancing decisions only happen in a few regions due to their computational complexity.}
\label{prediction}
\end{table}
\subsubsection{Performance of uncertainty prediction model}

We assessed various distributional assumptions regarding the uncertainty in demand. These probabilistic assumptions include Normal, Truncated Normal (TNorm), Poisson, Zero-Inflated Poisson (ZPoisson), and Negative Binomial (NB). These distributions span a diverse range of possible distributional assumptions, encompassing various categories such as continuous and discrete, real-line and non-negative. Their performances are summarized in \cref{prediction}.
\begin{table}[h]
\centering
\begin{tabular}{|c|c|c|c|c|c|}
\hline
\textbf{Distribution} & \centering\arraybackslash \textbf{NLL} & \centering\arraybackslash \textbf{MAE} & \centering\arraybackslash \textbf{MAPE} & \centering\arraybackslash \textbf{MPIW} & \centering\arraybackslash \textbf{PICP}\\ \hline
  \centering\arraybackslash \textbf{Normal} & 157 & \textbf{2.728} & \textbf{0.221} & 14.38 & 0.95\\ \hline
  \centering\arraybackslash \textbf{TNormal} & \textbf{153} & 2.745 & 0.223 & 15.13 & 0.88\\ \hline
  \centering\arraybackslash \textbf{Poisson} & 156 & 2.733 & 0.222 & 11.97 & 0.95\\ \hline
  \centering\arraybackslash \textbf{ZPoisson} & 156 & 2.737 & 0.222 & \textbf{11.92} & 0.95\\ \hline
  \centering\arraybackslash \textbf{NB} & 160 & 2.808 & 0.228 & 17.50 & \textbf{0.99}\\ \hline
\end{tabular}
\caption{Evaluation metrics for different probabilistic assumptions. The number cells in bold denote the best value in the column.}
\label{prediction}
\end{table}


\cref{prediction} reveals close composite (NLL) and point prediction (MAE and MAPE) accuracy among most of the tested distributions. However, distinctions emerge concerning MPIW and PICP. Notably, Poisson and Zero-Inflated Poisson exhibit the overall best fit with lower MPIW and higher PICP. It is noteworthy that the performance of Poisson and Zero-Inflated Possion is highly similar. This similarity arises from the fact that Zero-Inflated distributions are designed to accommodate datasets with a significant number of zeros. Given the limited presence of zeros in the dataset, the Zero-Inflated Poisson does not gain a significant advantage over the regular Poisson. Considering the dataset characteristics, Poisson is a more suitable choice, given its fewer parameters to estimate compared to Zero-Inflated Poisson. Therefore, we use the Poisson distribution across all DURO cases in the following experiments.
\subsubsection{Performance of \gls{abk:duro}}

\begin{figure}[tb]
\centering

\subfloat[Average waiting time of served passengers]{
	\includegraphics[width=0.5\textwidth]{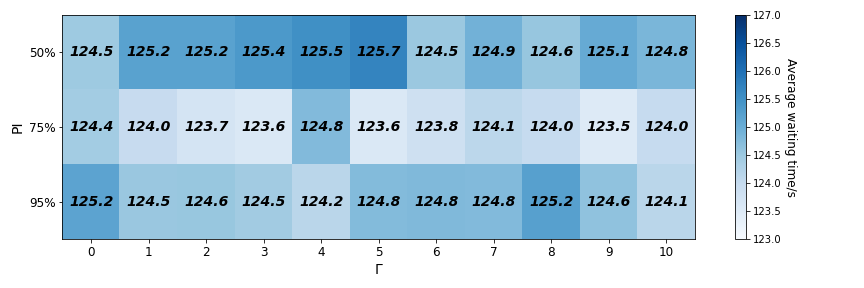} } 

\subfloat[Average traveling time of served passengers]{
	\includegraphics[width=0.5\textwidth]{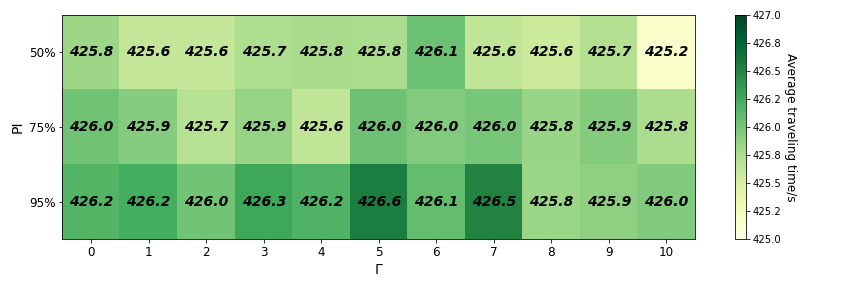} } 

\subfloat[Passenger leaving rate]{
	\includegraphics[width=0.5\textwidth]{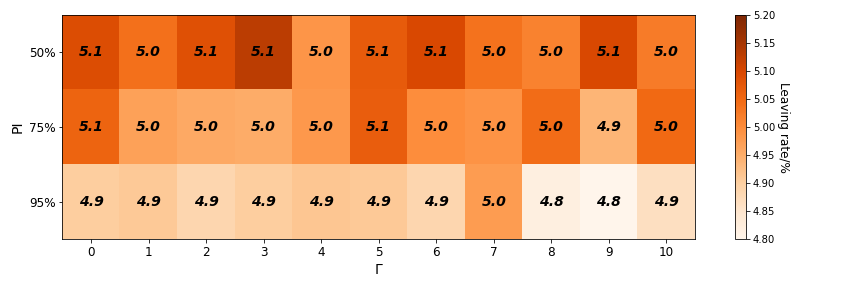}} 
	
\caption{Performance of DURO based on Poisson distribution. X-ais is the uncertainty parameter used in Constraint (\ref{eq:6b}), and the Y-axis is the percentage interval of uncertainty. These results serve as the baseline for later comparison.}
\label{fig: proposed}
\end{figure}

\cref{fig: proposed} shows the \gls{abk:duro} simulation result under the Poisson distribution assumption during the demand uncertainty quantification phase. The X-axis represents the value of the parameter $\Gamma$ in Constraint (\ref{eq:6b}) outlined in \cref{ro model}. It indicates the allowed fluctuation in the sum of the demand across all regions. The Y-axis corresponds to the Percentage Interval (PI) used for uncertain interval calculation. We experiment at 50\%, 75\%, and 95\% PI value.

The framework demonstrates a uniform performance across various experiments employing diverse parameter values. This observation implies that the neural networks consistently yield robust uncertainty estimation under varying PI values. The lowest average waiting time is attained with the 75\% interval. This finding indicates a trade-off between the coverage percentage and the interval width: elevating the PI to 95\% expands the encompassed range of potential demand patterns at the cost of introducing a wider interval that yields more conservative rebalancing decisions; reducing the PI to 50\% allows the optimization model to focus on the most probable scenarios, but at the potential expense of excluding pertinent cases. It is essential to underscore that the identified optimal PI value is not universally applicable to similar problems. Rather, the determination of the optimal PI is problem-specific and depends on the available data and the underlying probabilistic assumption.

The average traveling time of the served passengers is included as auxiliary information in the simulation process. It does not serve as a direct indicator of the quality of rebalancing decisions. Therefore, we focus solely on evaluating results based on the average waiting time and passenger leaving rate for the following comparative analysis. These two metrics more directly reflect the efficacy and impact of rebalancing decisions on service quality and passenger experience.
\subsubsection{Comparison with \gls{abk:dohv}} The historical average demand is a reasonable approximation of the true demand in cases with regular demand patterns. As illustrated in \cref{fig:demand} (b), the historical average closely aligns with the true demand in the dataset we employ. In \gls{abk:dohv}, we use the average demand before June 27 as an approximation of the true demand on June 27. This approximation serves as the input to the deterministic rebalancing model. 

\begin{figure}[tb]
\centering

\subfloat[Percentage of average waiting time reduction]{
	\includegraphics[width=0.5\textwidth]{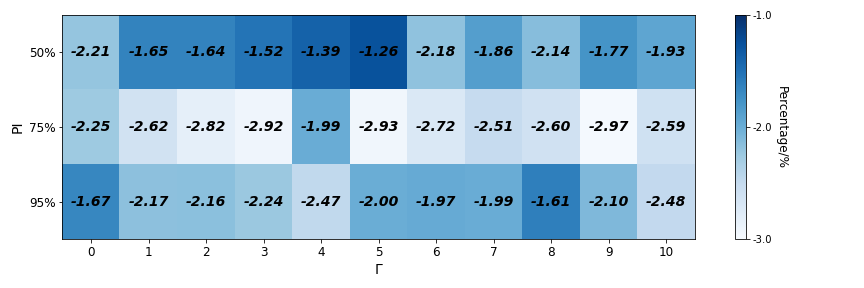} } 

\subfloat[Percentage of passenger leaving rate reduction]{
	\includegraphics[width=0.5\textwidth]{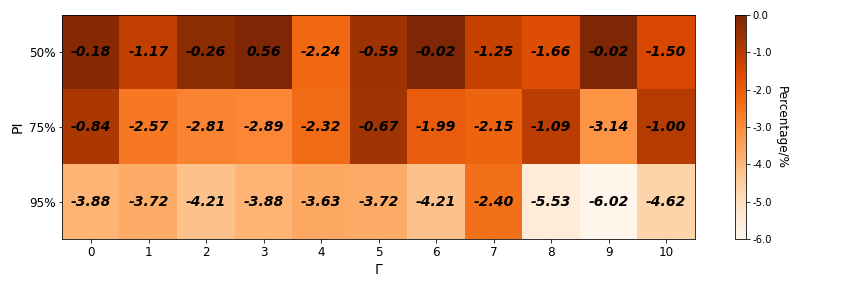}} 
	
\caption{Comparison between \gls{abk:duro} and \gls{abk:dohv}. \gls{abk:duro} achieves a lower average waiting time and passenger leaving rate across all instances compared to \gls{abk:dohv}.}
\label{fig: historical}
\end{figure}

\cref{fig: historical} shows the reduced percentage of the passenger waiting time and leaving rate achieved by \gls{abk:duro} in comparison to \gls{abk:dohv}. Across all simulated cases, \gls{abk:duro} consistently outperforms \gls{abk:dohv}, achieving a lower average waiting time and passenger leaving rate. This highlights the efficacy of the robust optimization model in delivering superior solutions when contrasted with its deterministic counterpart, even in scenarios where demand prediction is accurately presented.

\subsubsection{Comparison with \gls{abk:donn}} 
\begin{figure}[tb]
\centering

\subfloat[Percentage of average waiting time reduction]{
	\includegraphics[width=0.5\textwidth]{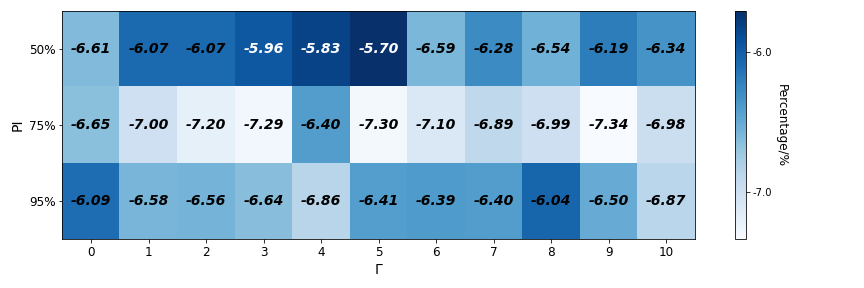} } 

\subfloat[Percentage of passenger leaving rate reduction]{
	\includegraphics[width=0.5\textwidth]{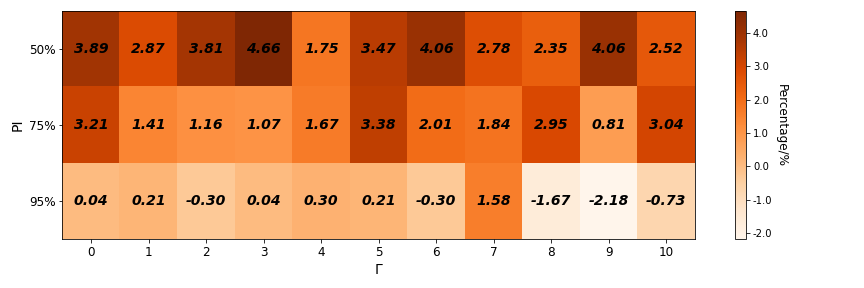}} 
	
\caption{Comparison between \gls{abk:duro} and DONN. Compared to DONN, \gls{abk:duro} reduces average waiting time across all cases, and lower leaving rate in the best-performing case.}
\label{fig: point}
\end{figure}

In contrast to the historical average, neural networks exhibit enhanced stability in prediction performance across diverse demand patterns and are widely adopted for demand forecasting. Demand prediction with deterministic optimization represents a leading-edge framework for addressing the vehicle rebalancing problem. To benchmark our framework against this approach, we trained and fine-tuned a graph-\gls{abk:lstm}-based neural network using demand data before June 27. This model was subsequently employed to forecast demand for June 27 within the 7:00 am to 9:00 am timeframe. The predicted demand was used as input to the deterministic counterpart of the model outlined in \cref{ro model}, excluding uncertainty constraints. The results are compared with those derived from \gls{abk:duro}, as displayed in \cref{fig: point}.

The average waiting time of \gls{abk:duro} outperforms DONN across all parameter combinations, as evidenced by a superior average waiting time. Despite DONN attaining a reduced passenger leaving rate in certain instances compared to \gls{abk:duro}, the best-performing case of \gls{abk:duro} still achieves a lower leaving rate. This reveals the superiority of \gls{abk:duro} over the deterministic optimization model, even when the latter also leverages the capabilities of neural networks.
\subsubsection{Comparison with RO}
With the RO rebalancing engine, the uncertainty set is defined by two parameters, $\Gamma$ and $\rho$. This subsection conducts a comparative analysis of RO's performance with various parameter combinations against \gls{abk:duro}. For \gls{abk:duro}, the result with $PI=75\%$ and $\Gamma=9$ is employed. This parameter combination achieves the lowest average waiting time, as identified from Figure \ref{fig: proposed}.

\begin{figure}[tb]
\centering

\subfloat[Percentage of average waiting time reduction]{
	\includegraphics[width=0.5\textwidth]{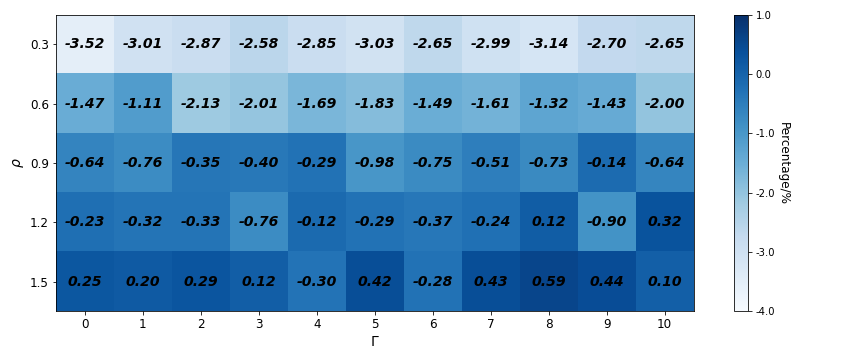} } 

\subfloat[Percentage of passenger leaving rate reduction]{
	\includegraphics[width=0.5\textwidth]{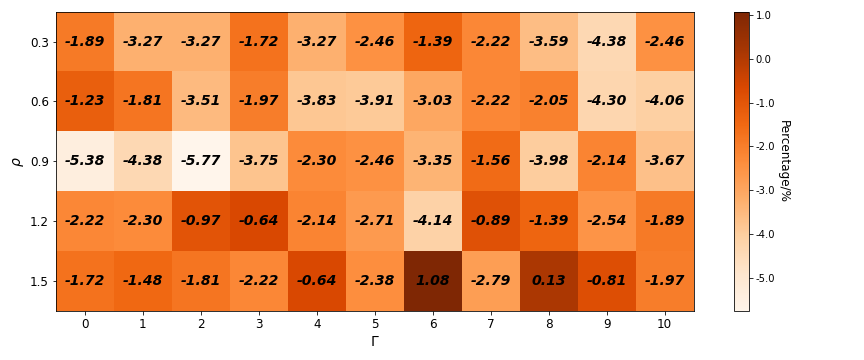}} 
	
\caption{Comparison between \gls{abk:duro} and RO. In most of the instances, \gls{abk:duro} achieves lower average waiting time and leaving rate compared to RO, and it eliminates the need for tuning $\rho$ in RO.}
\label{fig: ro}
\end{figure}

Figure \ref{fig: ro} shows that the best-performing model of \gls{abk:duro} achieves an overall lower average waiting time and leaving rate. Although in the case where $\Gamma = 8$ and $\rho = 1.5$, the RO model exhibits lower average waiting time and leaving rate, it is essential to underscore that \gls{abk:duro} eliminates the need for fine-tuning the values of $\rho$. In RO, $\rho$ is used to determine the uncertainty set, but the uncertainty quantification module in \gls{abk:duro} directly provides these intervals. Although the selection of the PI value is still necessary for \gls{abk:duro}, its meaning and interpretation are more clear and intuitive than $\rho$, making it easier to select a suitable PI value according to the application needs without an exhaustive iterative process. 

Additionally, it is pertinent to highlight that for the dataset used in the numerical experiments, the historical average demand already serves as a highly accurate approximation of the true demand, as demonstrated in Figure \ref{fig:demand}(b). Considering the robust generalization ability of neural networks, \gls{abk:duro}'s uncertainty quantification module can offer a better approximation of the true demand in cases where historical average demand falls short. This results in improved uncertainty sets and rebalancing optimization outcomes. In summary, the best-performing \gls{abk:duro} can achieve overall better results to \gls{abk:ro} with reduced tuning effort, and \gls{abk:duro} holds the potential for consistently superior performance in scenarios with irregular demand patterns.
\subsubsection{Comparison with \gls{abk:dro}}
\gls{abk:dro} has gained increasing popularity in optimization under uncertainty due to its flexibility in handling the uncertainty of the variables when the true distribution is not precisely known. \gls{abk:dro} eliminates the need for explicit distribution assumptions, allowing the model to infer the distribution from the provided data. In this subsection, we compare the performance of \gls{abk:duro} with that of \gls{abk:dro}. 

In the optimization module in \gls{abk:duro}, uncertainty is defined by possible value sets. In \gls{abk:dro}, the uncertainty is represented by sets of possible distributions. This introduces a higher level of complexity during the optimization process, which makes \gls{abk:dro} less computationally tractable. In our experiments, with the original problem size, we encounter memory constraints preventing the loading of the \gls{abk:dro} model.

To facilitate a meaningful comparison, we reduced the number of rebalancing regions from 63 to 30, focusing on high-demand areas centered around the south of Manhattan. Simulations were then executed with these reduced rebalancing regions for both \gls{abk:duro} and \gls{abk:dro} and the results were compared. Note that only the rebalancing operation is restricted within these 30 regions. The passenger trips and vehicle movements still cover the original 63 regions.

Constructing the ambiguity set from historical demand over the previous 18 days posed challenges for \gls{abk:dro}. The constructed ambiguity set results in an infeasible gls{abk:dro} problem because the historical data before June, 26th from 7:00 am to 9:00 am is insufficient to represent enough possible demands on June 26th during the same period. For example, if the demand in one region at 7:00 is zero for all the previous days, the lower and upper bounds of the demand on June, 26th at 7:00 for this region will be both zero. If the true demand for this region at the corresponding time is not zero, gls{abk:dro} will not be able to cover this possibility. Since gls{abk:dro} constructs the ambiguity set also based on the covariance of the uncertain variables, the insufficient quantification of the uncertainty set for one variable can be propagated to others, thus leading to an infeasible optimization problem. To address this, instead of using only the historical data at the same time for each region, we augmented the demand data with data from the previous and subsequent 5 minutes, assuming minimal demand changes during this interval. This enriched data is used to construct the ambiguity set encompassing a broader spectrum of possible demand scenarios. 

\begin{figure}[tb]
\centering

\subfloat[Percentage of average waiting time reduction]{
	\includegraphics[width=0.5\textwidth]{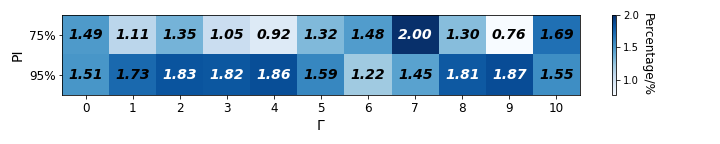} } 

\subfloat[Percentage of passenger leaving rate reduction]{
	\includegraphics[width=0.5\textwidth]{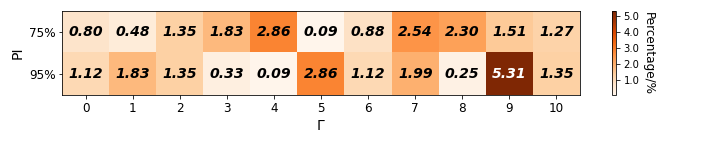}} 
	
\caption{Comparison between \gls{abk:duro} and \gls{abk:dro}. \gls{abk:dro} outperforms \gls{abk:duro} with a lower average waiting time and leaving rate.}
\label{fig: dro}
\end{figure}

The enriched dataset results in a solvable the \gls{abk:dro} model. The results are presented in \cref{fig: dro}. Due to changes in the rebalancing regions, the original 50\% interval from the uncertainty quantification module of \gls{abk:duro} is too conservative for the new problem, and the optimization module is not able to provide a solution based on these intervals. Consequently, experiments were conducted with 75\% and 95\% intervals.

\gls{abk:dro} attains a lower average waiting time and leaving rate compared to \gls{abk:duro}, attributed to its enhanced flexibility in handling uncertainty. \gls{abk:duro} makes probabilistic assumptions about the demand distribution, and uncertainty sets are derived based on these assumptions. This approach implicitly excludes distributions that are not tested but are potentially more consistent with the true distribution. In contrast, \gls{abk:dro} does not rely on a predefined distribution. It offers a robust solution against a wide range of possible uncertainty distributions. This feature enhances its flexibility in accommodating diverse uncertainty structures and makes it more robust when the true nature of uncertainty is not precisely known or is subject to change. Despite \gls{abk:dro}'s superior overall performance, \gls{abk:duro} can still achieve comparable results when the percentile and $\Gamma$ are appropriately selected. 

A notable disadvantage of \gls{abk:dro} is its significantly lower computation efficiency compared to \gls{abk:duro}. In our experiments, solving \gls{abk:dro} takes more than 60 times longer than \gls{abk:duro}. Moreover, as previously discussed, reducing the number of rebalancing regions in the experiments is necessary to make the \gls{abk:dro} model tractable with available computational memory. Even with the reduced problem size, \gls{abk:duro} takes less than 5s to solve for one time interval, whereas \gls{abk:dro} requires more than 5 minutes. Even accounting the parameters tuning time of \gls{abk:duro}, \gls{abk:dro} remains more time-consuming.

In summary, the analysis suggests that \gls{abk:dro} demands a larger and higher-quality dataset, and its computational efficiency is notably lower than \gls{abk:duro}. The choice between the \gls{abk:dro} and \gls{abk:duro} in real-world applications should depend on the available data and computation resources, with \gls{abk:duro} representing a more computationally efficient solution without compromising substantial performance gains.


\section{Conclusion}\label{conclusion}
The current research on optimization-based \gls{abk:amod} rebalancing often relies on complete information on future demand to derive solutions through deterministic optimization models. However, this approach frequently overlooks the inherent uncertainty in transportation systems. This paper introduces the \gls{abk:duro} framework, which is aimed for addressing this gap by leveraging a neural network to quantify the uncertainty about the knowledge on demand, and subsequently incorporating them into a robust optimization model. The framework is applied to the \gls{abk:amod} system rebalancing problem and showcases its competence in providing a reliable solution amid the demand uncertainty, as evidenced by comparisons with deterministic optimization, \gls{abk:ro}, and \gls{abk:dro} using real-world on-demand data from New York City.

The deterministic optimization models employed for comparison in our study adhere to a predict-and-optimize methodology. In this approach, the average demand is forecasted and input into a deterministic optimization model. Two distinct prediction methods are used: one based on historical averages and the other employing a graph-long short-term memory (graph-\gls{abk:lstm}) neural network. The outcomes of both models consistently demonstrate that \gls{abk:duro} achieves superior results with lower passenger waiting time and leaving rates. The comparison underscores the limitation of deterministic optimization, revealing its substantial dependence on point predictions, which restricts its ability to adapt to the inherent uncertainty in demand forecasting effectively.

In the comparison with \gls{abk:ro}, \gls{abk:duro} was compared with the uncertainty set determined by the historical mean and standard deviations, along with two manual selected parameters. Despite the best-performing \gls{abk:ro} achieving performance levels similar to \gls{abk:duro}, it is noteworthy that \gls{abk:ro} requires a more intensive effort in parameter tuning. While \gls{abk:ro} relies entirely on manual selected parameters to specify the uncertainty set, \gls{abk:duro} excels by utilizing neural network to offer a more efficient and streamlined process. 

While DRO exhibits higher flexibility in handling uncertainty, \gls{abk:duro} still attains comparable performance when the quantile and $\Gamma$ is appropriately selected. In addition, \gls{abk:duro} demonstrates higher computational efficiency and data generalization ability in comparison to DRO. This superior efficiency makes \gls{abk:duro} a more practical choice, especially when considering the constraints of available data and computational resources. This study underscores that \gls{abk:duro}, with its combination of performance and efficiency, stands out as a promising solution in tackling on-demand system rebalancing problem.

Future research directions include exploring the application of Bayesian quantile regression to enhance the uncertainty quantification module in \gls{abk:duro}. In the current framework, probabilistic assumptions are made on the demand distribution and uncertainty sets are calculated based on these assumptions. However, the reality of demand distribution may not conform perfectly to a distribution with well-defined mathematical characteristics. Bayesian quantile regression offers the potential to improve interval predictions, particularly in scenarios where demand distributions exhibit irregular patterns and resist easy characterization with traditional mathematical formulations. Furthermore, testing \gls{abk:duro} to on-demand systems characterized by higher uncertainty and irregular patterns is suggested to uncover its potential in more diverse scenarios. The demand data used in this paper exhibits high regularity across days. This regularity inadvertently made the momentum-based uncertainty set of \gls{abk:ro} and the ambiguity set specification in \gls{abk:dro} more straightforward, and obscure the full potential of \gls{abk:duro} when subjected to comparison.

\ifCLASSOPTIONcaptionsoff
  \newpage
\fi

\bibliographystyle{IEEEtran}
{\footnotesize
\bibliography{References.bib}}

\end{document}